\documentclass{elsart3-1}

\usepackage{amssymb}

\usepackage[english,francais]{babel}

\usepackage{graphicx}

\date{16 mars 2004}

\newfont{\bbf}{msbm10 at 12pt}
\newfont\bbfsm{msbm10 at 9pt}

\def\N{\mathbb {N}}

\def\Z{\mathbb {Z}}
\def\R{\mathbb {R}}
\def\C{\mathbb {C}}
\def\Cbar{\widehat{\C}} 
\def\Cstar{\C^*}

\def\disk{\mathbb {D}}

\def\diskstar{\disk^*}
\def\half{\mathbb{H}}

\newfont{\script}{eusm10 at 12pt}
\newfont{\scriptsmall}{eusm10 at 9pt}

\def\ovl{\overline}
\def\phi1{\phi}
\def\phi{\varphi}
\def\eps{\varepsilon}
\def\theta{\vartheta}
\def\Re{\mbox{\rm Re}}
\def\Im{\mbox{\rm Im}}

\def\sm{\setminus}

\def\K{\mathbb {K}} 
\def\k{ {\mbox{\tt k}} }
\def\*{ {\mbox{\tt $\star$}} }
\def\0{ {\mbox{\tt 0}} }
\def\1{ {\mbox{\tt 1}} }
\def\2{ {\mbox{\tt 2}} }
\def\3{ {\mbox{\tt 3}} }
\def\4{ {\mbox{\tt 4}} }

\def\s{ {\mbox{\tt s}} }

\def\St#1#2 {{\small\tt\rule{0pt}{0pt}_{\mbox{#1}}^{\mbox{#2}} }}

\def\<{\prec}
\def\>{\succ}

\ifx\ThmNum\undefined 
   \newtheorem{theorem}{Theorem}[section]
   \def\newsection#1{\setcounter{theorem}{0} \section{#1}}
\else 
   \newtheorem{theorem}{Theorem} 
   \def\newsection#1{\section{#1}}
\fi

\newtheorem{proposition}[theorem]{Proposition}
\newtheorem{lemma}[theorem]{Lemma}
\newtheorem{definition}[theorem]{Definition}

\newtheorem{corollary}[theorem]{Corollary}

\newtheorem{example}[theorem]{Example}

\def\proofof #1 {\par\medskip\noindent {\sc Proof of #1. }}

\def\sketchof #1 {\par\medskip\noindent {\sc Sketch of proof of #1. }}

\def\reminder #1 {{\sf #1}}
\def\hide #1 {}
\long\def\longhide #1 {}

\newcommand{\El}{E_{\lambda}}

\newcommand{\Eln}{\El^{\circ n}}
\newcommand{\Ek}{E_{\kappa}}

\newcommand{\Ekn}{\Ek^{\circ n}}

\renewcommand{\s}{{\underline s}}

\renewcommand{\j}{{\tt j}}

\newcommand{\Ray}{G_\s}
\newcommand{\Rayp}{G_{\s'}}
\newcommand{\Raypp}{G_{\s''}}
\newcommand{\IPR}{\Gamma_{W,h}}
\newcommand{\IPRp}{\Gamma_{W',h'}}
\newcommand{\RayOne}{G_{\s^{(1)}}}
\newcommand{\RayTwo}{G_{\s^{(2)}}}

\renewcommand{\theta}{\vartheta}
\newfont\Euler{eusm10 at 12pt}
\newcommand{\Sym}{\mbox{\Euler S}}

\renewcommand{\u}{{\tt u}}

\newcommand{\kk}{\underline{\k}}

\renewcommand{\St}[2]{{\rule{0pt}{0pt}_{#1}^{#2}}}

\def\LabelCaption#1#2{
  	\caption{{\sl #2}}
   \label{#1}
}

\def\addr{\mbox{\rm addr}}

\hide{
\newtheorem{theorem}{Theorem}[section]
\newtheorem{lemma}[theorem]{Lemma}
\newtheorem{e-proposition}[theorem]{Proposition}

\newtheorem{e-definition}[theorem]{Definition\rm}
\newtheorem{remark}{\it Remark\/}

\newtheorem{theoreme}{Th\'eor\`eme}[section]

\newtheorem{proposition}[theoreme]{Proposition}

\setcounter{equation}{0}
}

\def\og{\leavevmode\raise.3ex\hbox{$\scriptscriptstyle\langle\!\langle$~}}
\def\fg{\leavevmode\raise.3ex\hbox{~$\!\scriptscriptstyle\,\rangle\!\rangle$}}

\begin{document}

\begin{frontmatter}


\selectlanguage{english}
\title{Hyperbolic Components in Exponential
Parameter Space}

\vspace{-2.6cm}

\selectlanguage{francais}
\title{Composantes hyperboliques dans
l'espace des applications exponentielles}


\selectlanguage{english}
\author{Dierk Schleicher}
\ead{dierk@iu-bremen.de}

\address{School of Engineering and Science,
International University Bremen, Postfach 750 561, D-28725 Bremen}

\begin{abstract}
We discuss the space of complex exponential maps $\Ek\colon
z\mapsto e^{z}+\kappa$. We prove that every hyperbolic component
$W$ has connected boundary, and there is a conformal isomorphism
$\Phi_W\colon W\to\half^-$ which extends to a homeomorphism of
pairs $\Phi_W\colon(\ovl W,W)\to(\ovl\half^-,\half^-)$. This
solves a conjecture of Baker and Rippon, and of Eremenko and
Lyubich, in the affirmative. We also prove a second conjecture of
Eremenko and Lyubich. 
{\it To cite this article: Dierk Schleicher, C. R.
Acad. Sci. Paris, Ser. I 336 (2003).}

\vskip 0.5\baselineskip

\selectlanguage{francais}
\noindent{\bf R\'esum\'e}
\vskip 0.5\baselineskip
\noindent
Nous \'etudions l'espace des applications exponentielles complexes
$\Ek\colon z\mapsto e^{z}+\kappa$. Nous d\'emontrons que pour
chaque composante hyperbolique $W$, le bord $\partial W$ est
connexe, et qu'il y a un isomorphisme biholomorphe $\Phi_W\colon
W\to\half^-$ qui s'\'etend en un hom\'eomorphisme de paires
$\Phi_W\colon(\ovl W,W)\to(\ovl\half^-,\half^-)$.
Ceci \'etablit une conjecture de Baker et Rippon, et de Eremenko
et Lyubich. D'autre part, nous d\'emontrons une autre conjecture
de Eremenko et Lyubich. {\it Pour citer cet article~: Dierk
Schleicher, C. R. Acad. Sci. Paris, Ser. I 336 (2003).}

\end{abstract}
\end{frontmatter}

\selectlanguage{francais}
\section*{Version fran\c{c}aise abr\'eg\'ee}

Dans l'espace des applications exponentielles $\Ek\colon z\mapsto
e^z+\kappa$, chaque composante hyperbolique $W\subset\C$ est
simplement connexe, et il y a un isomorphisme conforme
$\Phi_W\colon W\to\half^-$ (le demi-plan gauche) tel que
l'application des multiplicateurs $\mu\colon W\to\diskstar$
se d\'ecompose comme $\mu=\exp\circ\Phi_W$. 
Il est assez facile de voir que $\Phi_W$ s'\'etend en une
application continue $\Phi_W\colon (\ovl
W,W)\to(\ovl\half^-,\half^-)$. 
Notre r\'esultat principal est que ceci est un hom\'eomorphisme.
Pour chaque $h\in\R$, nous consid\'erons le {\em rayon interne}
$\IPR\colon\R^-\to W$, $\IPR(t)=\Phi_W^{-1}(t+2\pi ih)$. On
d\'emontre facilement que pour chaque $h$, $\IPR(t)$ a une limite
dans $\Cbar$ (comme $t\nearrow 0$); la partie difficile est de
d\'emontrer que la limite est dans $\C$.

Pour les suites bourn\'ees $\s\in\Sym:=\Z^{\N}$, nous
introduisons les {\em rayons param\'etriques \`a l'adresse
externe $\s$:} ces sont des courbes 
$\Ray\colon(0,\infty)\to\C$ tel que pour
$\kappa=\Ray(t)$, l'orbite de la valeur singuli\`ere $\kappa$ sous
$\Ek$ converge vers $\infty$ (``l'orbite singuli\`ere
s'\'echappe''). Plus pr\'ecisement,
$\Ekn(\kappa)=F^{\circ n}(t)+2\pi i s_{n+1}+o(1)$, o\`u
$F(t)=e^t-1$ et $\s=s_1s_2s_3\dots$. En particulier, 
$\Ray(t)=t+2\pi i s_1+O(e^{-t})$. Les rayons param\'etriques ont
un ordre vertical naturel dans leur approche vers $+\infty$;
cet ordre est le m\^eme que l'ordre lexicographique de leurs
adresses externes $\s$.

Si, pour $W$ et $h$ donn\'ees, le rayon interne $\IPR$
aboutit \`a $\infty$, alors $\Re(\IPR(t))\to +\infty$, et
$\IPR$ d\'ecoupe l'espace des adresses externes en deux
ensembles $S^-$ et $S^+$ tels que $\Ray$ est dessus (ou dessous)
$\IPR$ si et seulement si $\s\in S^+$ (ou $\s\in S^-$). 
Pour la suite $\addr(\IPR):=\inf(S^+)=\sup(S^-)$, il y a
trois possibilit\'es: (1) $\addr(\IPR)\in\Sym$ est bourn\'ee;
(2) $\addr(\IPR)\in\Sym$ est non bourn\'ee; et 
(3) $\addr(\IPR)=s_1s_2\dots
s_{n-2}s_{n-1}$ est une suite finie telle que
$s_1,\dots,s_{n-2}\in\Z$ et
$s_{n-1}\in\Z+\frac 12$. Chacune de ces trois possibilit\'es
donnera une contradiction.

(1) Si $\s:=\addr(\IPR)$ est bourn\'ee, alors des calculs
asymptotiques impliquent qu'il n'y a pas de rayons
param\'etriques entre $\IPR$ et $\Ray$; mais l'ordre vertical et
des raisons combinatoires impliquent le contraire.

(2) Si $\addr(\IPR)$ n'est pas bourn\'ee, nous utilisons la
dynamique symbolique (en forme des {\em kneading sequences} et
des adresses internes) pour d\'emontrer qu'il y a une autre
composante hyperbolique $W'\neq W$ et deux rayons
param\'etriques $\RayOne$ et $\RayTwo$ qui aboutissent sur
$\partial W'$ et qui s\'eparent $W$ et $\IPR$ des rayons
param\'etriques aux adresses pr\`es de $\addr(\IPR)$, ce qui est
une autre contradiction.

(3) Si $\addr(\IPR)=s_1s_2\dots s_{n-2}s_{n-1}$ comme d\'ecrit
ci-dessus, alors il y a une composante hyperbolique $W'$ de
p\'eriode $n$ qui s'\'etend vers $\infty$, telle que le rayon 
param\'etrique $\Rayp$ s'approche vers $\infty$ dessus $W'$ ssi
$\s'\in S^+$. Encore une fois, il y a deux rayons 
param\'etriques $\RayOne$ et $\RayTwo$ qui aboutissent \`a $W'$ 
et qui s\'eparent $\IPR$ des rayons param\'etriques aux adresses
externes pr\`es de $\addr(\IPR)$, encore une contradiction.

Il s'en suit que $\IPR$ aboutit dans $\C$, et ceci suffit pour
d\'emontrer que $\Phi_W$ donne un hom\'eomorphisme de $\partial
W$ sur $\partial\half^-$, et que $\partial W$ est connexe. 

Finalement, nous d\'ecrivons deux autres conjectures de Eremenko
et Lyubich. Nous d\'emontrons qu'il existe une collection
d\'enombrable des composantes hyperboliques qui ne peuvent pas
\^etres jointes par des cha\^{\i}nes finies d'autres
composantes hyperboliques telles que les componsantes voisines
soient des bifurcations les unes des autres. Nous esp\'erons
que des m\'ethodes semblables \`a celles de notre d\'emonstration
du th\'eor\`eme~\ref{ThmBoundaryHypComps} pourraient aider \`a
d\'emontrer que des composantes non-hyperboliques sont bourn\'ees.

\selectlanguage{english}

\newsection {Introduction}
\label{SecIntro}

In this note, we investigate the fundamental structure of the
space of complex exponential maps $z\mapsto \Ek(z)=e^z+\kappa$
with $\kappa\in\C$. Translation by $-\kappa$ conjugates such a map
to $e^{z+\kappa}=\lambda e^z$ with $\lambda=e^\kappa$. The space
of complex exponential maps has been investigated since the
mid-1980's by Baker and Rippon \cite{BR}, Eremenko and Lyubich
\cite{ELR6T,ELP6,EL92}, Devaney, Goldberg and Hubbard
\cite{DGH}, and others.

A {\em hyperbolic component of period $n$} is a maximal open set
$W\subset\C$ such that for $\kappa\in W$, the map $\Ek$ has an
attracting periodic orbit of period $n$; all other periodic
orbits are then necessarily repelling. It is known from
\cite{EL92,BR,DGH} that every hyperbolic component is simply
connected, and it comes with a holomorphic multiplier map
$\mu\colon W\to\diskstar$ such that the attracting orbit of $\Ek$
has multiplier $\mu(\kappa)$. The map $\mu$ is a universal
covering map. Equivalently, there is a conformal isomorphism
$\Phi_W\colon W\to\half^-$ (the left half plane) such that
$\mu=\exp\circ\Phi_W$. The map $\Phi_W$ is unique up to
translation by $2\pi i\Z$. A preferred choice for $\Phi_W$ has
been given in \cite{AttrDyn}, but for our purposes any fixed
choice will do.

The main result of this note is the following.
\begin{theorem}
\label{ThmBoundaryHypComps} 
Every hyperbolic component has connected boundary,
and $\Phi_W$ extends to a homeomorphism of pairs $\mu\colon
(\ovl W,W)\to(\ovl\half^-,\half^-)$.
\end{theorem}
(Note that $\ovl W$ and $\ovl\half^-$ etc.\ will denote closures
in $\C$ throughout this paper; however, this theorem remains true
if closures in the Riemann sphere $\Cbar$ are taken.)
This result had been conjectured by Baker and Rippon~\cite{BR}
and Eremenko and Lyubich~\cite{ELR6T} in the mid-1980's.
The proof requires a substantial amount of knowledge on
exponential parameter space; many of the required results are of
interest in their own right. Among them is a description of
parameters $\kappa$ for which the singular orbit {\em escapes},
i.e.\ converges to $\infty$.
We need the map $F\colon\R^+\to\R^+$, $F(t)=e^t-1$ and the
notation $\Sym:=\Z^{\N}$; sequences $\s\in\Sym$ will be called
{\em external addresses}.

\begin{theorem}
\label{ThmBoundedParaRays}
For every bounded external address $\s\in\Sym$, there exists a
unique injective $C^1$-curve $\Ray\colon(0,\infty)\to\Cstar$ 
(a {\em parameter ray}) so that for $\kappa=\Ray(t)$, the
singular value $\kappa$ escapes to $\infty$ such that
\[
\Eln(\kappa)=F^{\circ n}(t)+2\pi is_{n+1}+o(1)
\qquad\mbox{as $n\to\infty$ or $t\to\infty$} \,\,.
\]
The curve $\Ray$
satisfies $\Ray(t)=t+2\pi i s_1+O(e^{-t})$ as $t\to\infty$.
All these curves are disjoint. 
\end{theorem}
In fact, there are parameter rays $\Ray$ for
all {\em exponen\-tially bounded} external addresses $\s\in\Sym$
\cite{Markus}. Exponential maps with escaping singular orbits
are completely classified in terms of parameter rays \cite{FRS}.

Note that all parameter rays $\Ray$ come with a natural vertical
order: since these rays are disjoint and $\Re(\Ray(t))=+\infty$
as $t\to+\infty$, each ray cuts sufficiently far right half
planes into two unbounded parts, so every other parameter ray
must be {\em above} or {\em below} $\Ray$ (depending on in which
unbounded part it converges to $+\infty$). The proof of
Theorem~\ref{ThmBoundedParaRays} also shows that the vertical
order coincides with the lexicographic order of the external
address $\s$. We say that {\em the parameter ray $\Ray$ lands at
$\kappa\in\C$} if $\lim_{t\searrow 0}\Ray(t)=\kappa$.

A periodic orbit is {\em indifferent} if it has a periodic orbit
with multiplier $\mu\in\partial\disk$. Every such parameter is on
the boundary of a hyperbolic component. If $\mu$ is a root of
unity, the parameter is called {\em parabolic}.

\begin{theorem}
\label{ThmPeriodicParaRays}
For every periodic $\s\in\Sym$, the parameter ray $\Ray$ lands at
a parabolic parameter $\kappa$, and every parabolic $\kappa$ is
the landing point of one or two parameter rays at periodic
external addresses.
\end{theorem}
For the purposes of this note, only the second half of
Theorem~\ref{ThmPeriodicParaRays} will be needed; the proof of
the first half needs Theorem~\ref{ThmBoundaryHypComps}.
(Similarly, if $\s$ is preperiodic, then $\Ray$ lands at a {\em
Misiurewicz parameter:} that is a $\kappa$ such that the singular
orbit of $\Ek$ is strictly preperiodic; conversely, every
Misiurewicz parameter is the landing point of a finite positive
number of parameter rays at preperiodic external addresses.)

Most results in this note were obtained in the thesis
\cite{Habil}. A detailed proof of
Theorem~\ref{ThmBoundaryHypComps} will be given in
\cite{LRDS}, together with the necessary background about
exponential parameter space and a number of further results.

\section{Internal Rays}

Consider a hyperbolic component $W$ with conformal isomorphism
$\Phi_W\colon W\to\half^-$ as above. For every $k\in\Z$, the set
$\Phi^{-1}(\{z\in\half^-\colon 2\pi k<\Im(z)<2\pi(k+1)\}$ is
called a {\em sector} of $W$. 
Moreover, for any $h\in\R$ we define the {\em internal ray at
height $h$} to be the curve $\IPR\colon\R^-\to W, t\mapsto 
\Phi_W^{-1}(t+2\pi ih)$. The component
$W$ is thus canonically foliated into internal rays.

\begin{lemma}
The multiplier map $\mu\colon W\to\diskstar$ and
the conformal isomorphism $\Phi_W\colon W\to\half^-$
extend to continuous maps $\mu\colon\ovl{W}\to\ovl{\diskstar}$ and
$\Phi_W\colon\ovl W\to \ovl\half^-$.
Every boundary component of $W$ is a piecewise analytic curve.
For every $W$ and $h$, $\lim_{t\to-\infty}\IPR(t)=+\infty$,
while $\lim_{t\nearrow 0}\IPR(t)$ exists in $\Cbar$.
\end{lemma}
The first two statements follow simply from the implicit
function theorem. Since as $t\to-\infty$, $\mu(\IPR(t))\to 0$,
but no exponential map has a periodic orbit with multiplier $0$,
it follows
$\lim_{t\to-\infty}\IPR(t)=\infty$. As $t\nearrow 0$, any limit
parameter in $\C$ of $\IPR(t)$ must have an indifferent orbit
with multiplier $\exp(h)$; since such are discrete, the ray
{lands} in $\Cbar$. Theorem~\ref{ThmBoundaryHypComps} is proved
as soon as we know that all limits are in $\C$: this gives a
continuous inverse $\Phi_W^{-1}\colon\ovl\half^-\to\ovl W$. 

In order to prove that the internal ray $\IPR$ lands in
$\C$ (as $t\nearrow 0$), we use the (external) parameter rays
$\Ray$ from Theorem~\ref{ThmBoundedParaRays}. Suppose that $\IPR$
lands at $\infty$ as $t\nearrow 0$. It is not hard to show
that the real parts tend to $+\infty$ (the ray cannot cross
other hyperbolic components or parameter rays $\Ray$), so
the vertical order between the ray $\IPR$ and each parameter ray
$\Ray$ is well-defined. Therefore, $\IPR$ cuts the space
of bounded external addresses into two sets $S^+$ and $S^-$ such
that rays $\Ray$ with $s\in S^+$ are above $\IPR$, while rays with
$s\in S^-$ are below $\IPR$. This defines a cutting
sequence $\addr(\IPR):=\inf\{S^+\}=\sup\{S^-\}$ for which
there are three possibilities:
\begin{description}
\item[(1)]
$\addr(\IPR)$ is a bounded sequence in $\Sym$;
\item[(2)]
$\addr(\IPR)$ is an unbounded sequence in $\Sym$;
\item[(3)]
$\addr(\IPR)=s_1s_2\dots s_{n-2}s_{n-1}$, where
$s_1,\dots,s_{n-2}\in\Z$, while $s_{n-1}\in\Z+\frac 1 2$.
\end{description}
The third case needs some explanation: since $\addr(\IPR)$ is defined
as a supremum over bounded sequences in $S^-$, there might be a
position $n$ such that the supremum of the first $n-1$ entries is
a finite sequence $s_1s_2\dots s_{n-2}s'_{n-1}$ of integers,
while the $n$-th entries are unbounded above. Setting
$s_{n-1}:=s'_{n-1}+\frac 1 2$, the finite sequence $s_1s_2\dots
s_{n-2}s_{n-1}$ can then be considered the supremum over $S^-$,
as well as the infimum over $S^+$. In order to prove that the
internal ray $\IPR$ cannot land at $\infty$, we have to
exclude all three possibilities for $\addr(\IPR)$.

\section{Squeezing of Internal Rays}

Given an internal ray $\IPR$ which lands at $\infty$, we exclude
the three cases for $\addr(\IPR)$ in order.

{\bf(1)} If $\s:=\addr(\IPR)$ is bounded, then there is a
parameter ray $\Ray\colon(0,\infty)\to\C$. To fix ideas, suppose
that $\IPR$ approaches $+\infty$ below $\Ray$. All parameter rays
$\Rayp$ with $\s'\in S^-$ are below $\IPR$. However, for every
$\eps>0$ and $s\in\N$ there is an $n\in\N$ such that
$|\Ray(t)-\Rayp(t)|<\eps$ uniformly for all $t>1$ provided the
first $n$ entries in $\s$ and $\s'$ coincide and all entries in
$\s$ and $\s'$ are bounded by $s$. No matter how closely $\IPR$
approaches $+\infty$ to $\Ray$, there is another ray $\Rayp$
closer to $\Ray$, and this is a contradiction.

{\bf (2)} If $\addr(\IPR)$ is an unbounded sequence in $\Sym$, we
need to know quite a bit more about exponential parameter space.
The fundamental idea is easy, though: 
there are a hyperbolic component $W'$ and two external addresses
$\s^{(1)}$ and $\s^{(2)}$ such that the parameter rays $\RayOne$
and $\RayTwo$ land at $\partial W'$, and $\IPR$ is in a different
connected component of $\C\sm(\ovl W'\cup\RayOne\cup\RayTwo)$
than all parameter rays $\Rayp$ at external addresses near
$\addr(\IPR)$.
This is a contradiction again.

The basic idea is to use symbolic dynamics in the form of {\em
kneading sequences} (and the human-readable variant, {\em
internal addresses}) \cite{IntAddr}. For fixed $\s\in\Sym$ and all
$\k\in\Z$, consider the sequences $\k\s$ (concatenation: the
symbol $\k$ followed by the sequence $\s$). Then
$\Sym\sm\bigcup_{\k}\{\k\s\}$ is the union of the countably many
intervals
$\Sym_{\k}=(\k\s,(\k+1)\s)$. We define the {\em kneading sequence}
$\K(\s)$ as the infinite sequence $\k_1\k_2\k_3\dots$ such that
$\k_i=\k$ iff $\sigma^{i-1}(\s)\in \Sym_{\k}$; in the boundary
case when $\sigma^{i-1}(\s)=\k\s$, set $\k_i=\St{\k-1}{\k}$.
Clearly, a boundary symbol $\St{\k-1}{\k}$ occurs if and only if
$\s$ is periodic, and $\K(\s)$ is bounded if and only if $\s$ is
bounded.

Here is one aspect how kneading sequences help to describe the
structure of parameter space.
\begin{proposition}
(a)
Every sector $W'_k$ of a hyperbolic component $W'$ of period
$n$ has an associated sequence $\kk=\k_1\k_2\k_3\dots$ which is
periodic of period $n$ with the following property: if the
parameter ray $\Rayp$ lands on $\partial W'_k$, and $\s'$ is
periodic, then $\K(\s')$ is ``almost
equal'' to $\kk$ in the following sense: the period of $\s'$
equals $qn$ for some $q\in\N$, and the $i$-th entry of $\K(\s')$
equals $k_i$ whenever $i$ is not a multiple of $qn$; if it is,
then the $i$-th entry of $\K(\s)$ is either
$\St{\k_i-1}{\k_i}$ or $\St{\k_i}{\k_i+1}$.

(b)
Suppose that $\s'$ and $\s''$ are two bounded external addresses
whose kneading sequences coincide in their first $n-1$ entries,
while the $n$-th entries differ. Then there there are a hyperbolic
component $W'$ of some period $n'\leq n$ and two external
addresses $\s^{(1)}$ and $\s^{(2)}$ such that the parameter rays
$\RayOne$ and $\RayTwo$ land at $\partial W'$, and $\Rayp$ and
$\Raypp$ are in different connected components of $\C\sm(\ovl
W'\cup\RayOne\cup\RayTwo)$.
\end{proposition}

We can now finish the proof that the internal ray $\IPR\subset W$
cannot land at $+\infty$ so that $\addr(\IPR)$ is unbounded. The
ray $\IPR$ is contained in the closure of some sector $W_k$ of
$W$, and all parameter rays landing at $\partial W_k$ have
uniformly bounded kneading sequences. Since $\addr(\IPR)$ is
unbounded, so is its kneading sequence, and there must be
infinitely many pairs of hyperbolic components with associated
parameter rays landing at them which separate $\IPR$ from all
parameter rays landing at $\partial W_k$.

{\bf (3) }
The third case is that $\addr(\IPR)=s_1s_2\dots
s_{n-2}s_{n-1}=:\s$, where $s_1,\dots,s_{n-2}\in\Z$, while
$s_{n-1}\in\Z+\frac 1 2$. By the main result of \cite{AttrDyn},
there exists a unique hyperbolic component $W'$ of period $n$
such that for every $h\in\R$, a parameter ray $\Rayp$
approaches $+\infty$ above (resp.\ below) the internal
ray $\IPRp(t)$ if and only if $\s'>\s$ (resp.\ $\s'<\s$) (recall
that $\IPRp$ satisfies $\lim_{t\to-\infty}\Re(\IPRp(t))=+\infty$).
This means that the ray $\IPR$ approaches $+\infty$ (as
$t\to 0$) so close to the curve $\IPRp$ (as $t\to-\infty$) that no
parameter ray at bounded external address is between $\IPR$ and
$\IPRp$. This is excluded by the following result.

\begin{lemma}
For every hyperbolic component $W'$ of period $n$ with associated
external address $s_1s_2\dots s_{n-2}s_{n-1}$, and for every
$\xi>0$, there are two parameter rays $\RayOne$ and $\RayTwo$
which both land at $\partial W'$, and the two landing points can
be connected by a curve $\Gamma\subset W$ such that all points in
$\RayOne\cup\RayTwo\cup\ovl\Gamma$ have real parts greater than
$\xi$.

The curve $\RayOne\cup\RayTwo\cup\ovl\Gamma$ disconnects $\C$ into
two parts, say $U$ and $U'$, so that all real parts in $U'$ are
greater than $\xi$, and all parameter rays with external addresses
between $\s^{(1)}$ and $\s^{(2)}$ are contained in $U'$.
\end{lemma}

If $W\neq W'$, then it is easy to see that $\IPR$ cannot approach
$+\infty$ so that $\addr(\IPR)=s_1s_2\dots s_{n-2}s_{n-1}$: the ray
$\IPR$ is disjoint from $\RayOne\cup\RayTwo\cup\ovl\Gamma$, so if
$\xi$ is sufficiently large, then $\IPR$ cannot be contained in
$U'$; but parameter rays at external addresses close to $\s$ must
be contained in $U'$, a contradiction.

Finally, if $W'=W$, then $W$ has associated external address 
$s_1s_2\dots s_{n-2}s_{n-1}$. The two ends of the curve $\IPR$
(as $t\to 0$ and $t\to-\infty$) are not homotopic within $W$, so
they enclose together some part of $\partial W$. But the
surrounded part of $\partial W$ contains infinitely many
parabolic parameters, hence infinitely many parameter rays at
periodic external addresses, and again $\addr(\IPR)\neq
s_1s_2\dots s_{n-2}s_{n-1}$.

This concludes the proof that every internal ray lands in $\C$,
and hence the proof of Theorem~\ref{ThmBoundaryHypComps}.

In addition to connectedness of the boundaries of hyperbolic
components, \cite{ELR6T} contains two more conjectures about
exponential parameter space. The first states that there are
countably many hyperbolic components of which no two can be
connected by a finite chain of further hyperbolic components so
that adjacent components in the chain are bifurcations from each
other. This conjecture also follows from a systematic
investigation of the bifurcation structure of hyperbolic
components as given in \cite{Habil,LRDS}; in fact, there are
countably many components of period $3$ with imaginary parts in
$(-\pi,+\pi)$, and among them is a single pair of
components which can be connected by a chain of bifurcating
components. 

The third conjecture in \cite{ELR6T} states that {\em
hyperbolicity is dense in exponential parameter space}, in
analogy to the main conjecture about quadratic polynomials. 
If this is false,
then there is a {\em non-hyperbolic component} $W\subset\C$: this
is a maximal open set containing no exponential map with an
attracting periodic orbit. A similar argument as above shows
that $W$ cannot contain a curve to $\infty$ \cite{LRDS};
moreover, there are at most two external addresses $\s_W$ and
$\s'_W$ such that every parameter ray $\Ray$ at external address
$\s\notin\{\s_W,\s'_W\}$ is separated from $W$ by a pair of
periodic parameter rays landing at a common point. 
It might well be possible to close these routes to $\infty$ for
$W$ in a similar way, thus proving that every non-hyperbolic
component (if any) was bounded. This nicely contrasts with the
observation that many features of exponential parameter space are
unbounded (such as all hyperbolic components). If one could prove
that every non-hyperbolic component had to be unbounded,
then this would prove the third conjecture.




\begin{figure}[h]
\includegraphics[height=80mm]{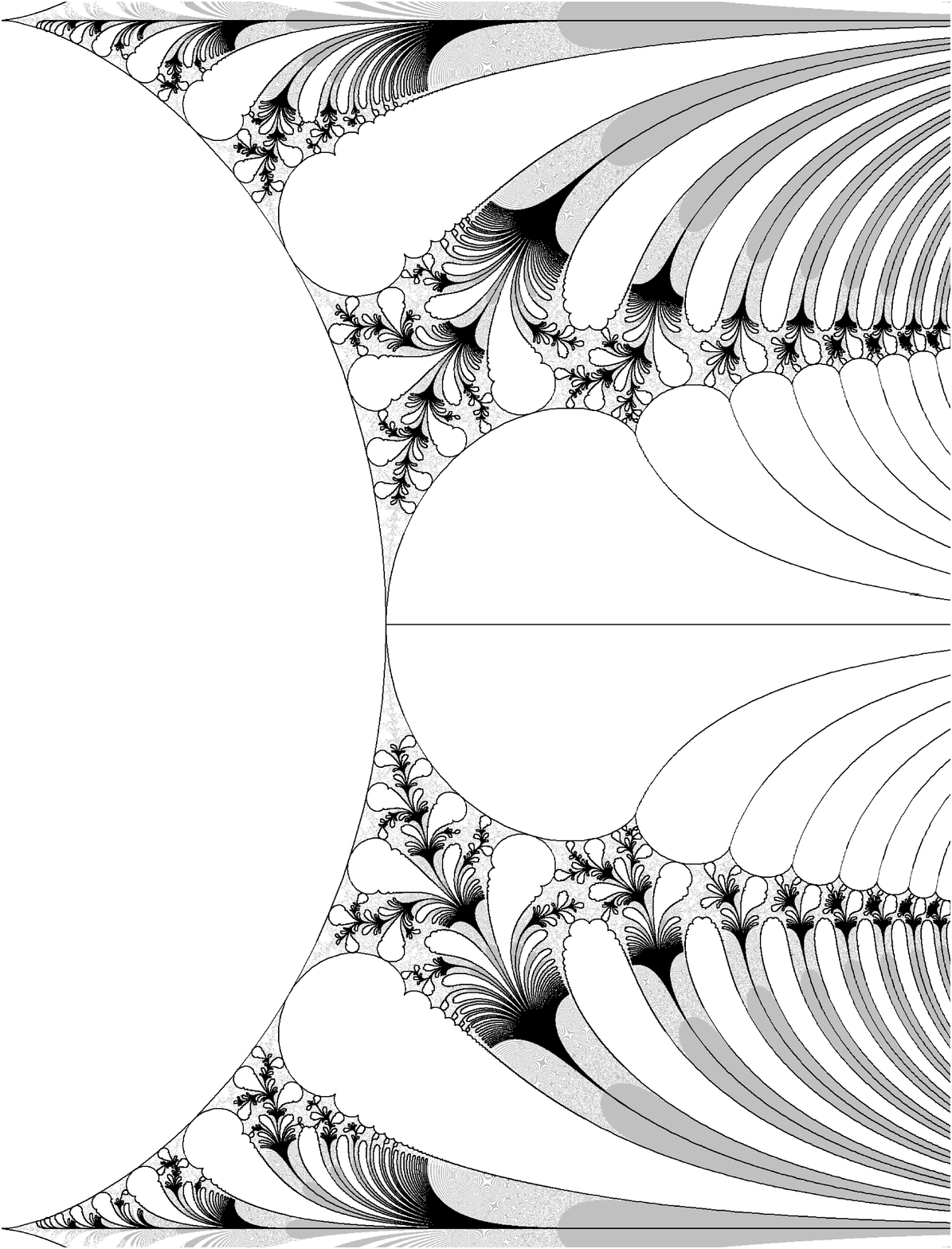}
\LabelCaption{FigExpoParaSpace}
{The space of complex exponential maps $z\mapsto e^z+\kappa$.
Many hyperbolic components of various periods are shown. It is
clearly visible that each has connected boundary. ---
L'espace des applications exponentials complexes $z\mapsto
e^z+\kappa$. Beaucoup des composantes hyperboliques sont
dessin\'ees. Il est bien visible que leur bord est connexe.}

\end{figure}


\begin{thebibliography}{MM}

\input{cyracc.def}

\def\j{{\u i}}
\def\J{{\u I}}
\newfont{\cyrit}{wncyi10 at 8pt}

\bibitem[BR]{BR}
I.~Noel Baker and Phil J.~Rippon:
{\em Iteration of exponential functions}. Ann.~Acad.\ Sci.\ Fenn.,
Series A.I.\ Math.\ {\bf 9} (1984), 49--77.

\bibitem[DGH]{DGH}
Robert Devaney, Lisa Goldberg and John Hubbard: {\em A dynamical
approximation to the exponential map by polynomials}.
Preprint, MSRI Berkeley (1986).

\bibitem[EL1]{ELR6T}
Alexandre Eremenko and Mikhail Lyubich:
{\em Iterates of Entire Functions}.
{ Soviet Math. Dokl.} {\bf 30} (1984), 592--594.

\bibitem[EL2]{ELP6}
Alexandre Eremenko and Mikhail Lyubich:
{{\cyracc \cyrit Iteratsii tselykh funktsii}}
{\em (Iterates of Entire Functions)}.
Preprint, Physico-Technical Institute of Low-Temperatures
Kharkov {\bf 6} (1984).

\bibitem[EL3]{EL92}
Alexandre Eremenko and Mikhail Lyubich:
{\em Dynamical properties of some classes of entire functions}.
Annales de l'Institut Fourier, Grenoble {\bf 42} 4 (1992),
989--1020.

\bibitem[FS]{Markus}
Markus F\"orster and Dierk Schleicher:
{\em Parameter rays for the exponential family}. In preparation.

\bibitem[FRS]{FRS}
Markus F\"orster, Lasse Rempe, and Dierk Schleicher:
{\em Classification of escaping exponential maps}. Submitted.

\bibitem[LS]{IntAddr}
Eike Lau and Dierk Schleicher: {\em Internal addresses in the
Mandelbrot set and irreducibility of polynomials}. Preprint {\bf
14}, Institute of Mathematical Sciences, Stony Brook (1994).

\bibitem[RS]{LRDS}
Lasse Rempe and Dierk Schleicher:
{\em Bifurcations in the space of exponential maps}. 
In preparation.

\bibitem[S1]{Habil}
Dierk Schleicher: {\em On the dynamics of iterated exponential
maps}. Habilitationsschrift, Technische Universit\"at M\"unchen
(1999).

\bibitem[S2]{AttrDyn}
Dierk Schleicher: {\em Attracting Dynamics of Exponential Maps}.
Ann.~Acad.\ Sci.\ Fenn., Series A.I.\ Math.\ {\bf 28} 1
(2003), 3--34.






\end{thebibliography}
\end{document}